\documentclass[a4paper, 12pt]{amsproc}
\usepackage{amssymb,amsthm,amsfonts,latexsym}
\theoremstyle{remark}

\title[Iwahori's Question for Affine Hecke Algebras]
{Iwahori's Question for Affine Hecke Algebras}

\author[T. Shoji and N. Xi]{Toshiaki Shoji$^{*}$ and Nanhua Xi$^{\dagger}$}
\address{$^{*}$
 Graduate School of Mathematics\\
Nagoya University\\
Chikusa-ku, Nagoya\\
464-8602 Japan} \email{shoji@math.nagoya-u.ac.jp }
\address{$^{\dagger}$
Institute of Mathematics\\
Chinese Academy of Sciences\\
Beijing, 100190\\
China } \email{nanhua@math.ac.cn}
\begin{document}
\baselineskip=18pt
\begin{abstract}
In this paper we show that  an affine Hecke algebra $H_q$ over
complex numbers field  with parameter $q\ne 1$ is not isomorphic to
the group algebra over complex numbers field of the corresponding
extend affine Weyl group if the corresponding root system has no
factors of type $A_1$ and the order of $q$ is different from
11 and 13 if the root system has factors of type $E_8$.

\end{abstract}

\maketitle

\def\Cal{\mathcal}
\def\bold{\mathbf}
\def\ca{\mathcal A}
\def\cdz{\mathcal D_0}
\def\cd{\mathcal D}
\def\cdo{\mathcal D_1}
\def\bold{\mathbf}
\def\l{\lambda}
\def\le{\leq}

\def\ds{\displaystyle\sum}
\def\vp{\varphi}
\def\st{\stackrel}
\def\sc{\scriptstyle}

\def\bbQ{\mathbb Q}
\def\bbF{\mathbb F}
\def\bbC{\mathbb C}

\def\bbN{\mathbb N}

For a Hecke algebra (over a field) of a finite Coxeter group, it is
known that the Hecke algebra is isomorphic to the group algebra over
the field of complex numbers when both the Hecke algebra and the
group algebra are semisimple. For a simple criterion of the
semi-simplicity for the Hecke algebra we refer to [G]. It is natural
to consider the question for affine Hecke algebra, this question is
analogue to a question of Iwahori for Weyl group[I]. For affine
Hecke algebras the answer is different. In an unpublished work
around 1970's, Casselman proved, by using group cohomology, that the
Hecke algebra of a simple group over a $p$-adic field $k$ with
respect to an Iwahori subgroup is not isomorphic to the group
algebra of the corresponding affine Weyl group. unless the group is
of $k$-rank one. (The authors thank the referee for pointing out
this fact.) We will show that an affine Hecke algebra $H_q$ over the
complex number field with parameter $q\ne 1$ is not isomorphic to
the group algebra over complex numbers field of the corresponding
extend affine Weyl group if the the root system has no factors of
type $A_1$. For technical reason, we require that the order of $q$
is different from 11 and 13 if the root system has factors of type
$E_8$, see Theorem 1.2.

\section{Affine Hecke algebras}

\noindent {\bf 1.1.} Let $R$ be an irreducible root system, $W_0$
the Weyl group of $R$, $Q=\mathbb Z R$ the root lattice and $X$ the
weight lattice. The Weyl group $W_0$ acts on $Q$ and $X$. Then the
semidirect product $W_a=W_0\ltimes Q$ is an affine Weyl group, which
is a subgroup of the extend affine Weyl group $W=W_0\ltimes X$. Fix
a positive root system of $R$ and denote by $R^+$ the set of
positive roots of $R$. Then we have a length function $l:
W\to\mathbb N$ given  by the formula (see \cite{IM})
$$l(wx)=\displaystyle\sum_{\st{\sc \alpha\in
R^+} {w(\alpha)\in R^-}}|\langle
x,\alpha^\vee\rangle+1|+\displaystyle\sum_{\st{\sc \alpha\in R^+}
{w(\alpha)\in R^+}}|\langle x,\alpha^\vee\rangle|,$$ where
$R^-=-R^+$. The set $S$ of simple reflections consists of all
elements in $W_a$ with length 1. The set of dominant weights $X^+$
is $\{x\in X\ |\ l(wx)=l(w)+l(x)\}$.

\def\Cal{\mathcal}
\def\bold{\mathbf}
\def\ca{\mathcal A}
\def\cdz{\mathcal D_0}
\def\cd{\mathcal D}
\def\cdo{\mathcal D_1}
\def\bold{\mathbf}
\def\l{\lambda}
\def\ca{\Cal A}
\def\cb{\Cal B}
\def\hkq{H_{q}}
\def\hq{H_q}
\def\zhkq{Z(\hkq)}
\def\Om{\Omega}
\def\om{\omega}
\def\End{\text {End}}
\def\vp{\varphi}
\def\st{\stackrel}
\def\sc{\scriptstyle}
\def\bbC{\mathbb C}

Denote by $\hkq$ the Hecke algebra of $(W,S)$ over the complex
numbers field $\bbC$ with a non zero parameter $q\in \bbC^*$. By
definition, $\hkq$ has a $\bbC$-basis consisting of elements $T_w,\
w\in W$, and the multiplication law is given by the relations:
$(T_r-q)(T_r+1)=0$ if $r\in S$, $T_wT_u=T_{wu}$ if
$l(wu)=l(w)+l(u)$.

The main result of this article is the following.

\medskip

 \noindent {\bf Theorem 1.2.} Let $q\ne 1$ be a nonzero complex number and $R$ an irreducible
 root system of rank greater than 1. Assume that the order of $q$ is different
from 11 and 13 if $R$ is of type $E_8$.  Then $\hkq$ is not
isomorphic to $\bbC[W]$.

\medskip

{\it Remark.} When $W$ is of type $\tilde A_2$, the result was
proved in [X1, 11.7]. When $W$ is of type $\tilde A_1$, we know that
$\hkq$ is isomorphic to $\bbC[W]$ if and only if $q\ne -1$, see
loc.cit.. For type $\tilde A_2$,  Yan showed that $\hkq$ is
isomorphic to $H_p$ if and only if $p=q$ or $pq=1$.

The theorem should be valid for $R$ being of type $E_8$ and the
order of $q$ is 11 or 13, but the authors have not been able to work
out it. Recently the authors learned that this result is known to G.
Lusztig for long time.

\medskip

\noindent{\bf 1.3.} The center of $\hkq$ will play a key role in the
proof. We need to recall the description of Bernstein for the
center.

For each $x$ in $X$, we can find $y$ and $z$ in $X^+$ such that
$x=yz^{-1}$. Define $\theta_x=q^{\frac12(l(z)-l(y))}T_yT_z^{-1}$. It
is known that $\theta_x$ is independent of the choice of $y$ and
$z$.

\def\Tq{\Theta_q}

\medskip

\noindent (a) (Bernstein) For any $x,y$ in $X$,
$\theta_x\theta_y=\theta_y\theta_x=\theta_{xy}$. The elements
$\theta_x, \ x\in X$ form a $\bbC$-basis of the subalgebra
$\Theta_q$
 of $\hkq$ generated by all $\theta_y, \ y\in X$. The Hecke algebra $\hkq$ is a free $\Tq$-module with a basis $T_w,\ w\in W_0$.

\medskip

  For each dominant weight $x$ in
$X$, let $O_x$ be the $W_0$-orbit of $x$. That is
$$O_x=\{wxw^{-1}\in X\ |\
 w\in W_0\}.$$ For $x\in X^+$, define
$$S_x=\sum_{y\in O_x}\theta_y.$$ Let $n$ be the rank of $R$ and
$x_1,x_2,...,x_n$ be the fundamental dominant weights. Then we have
(see [L1, Theorem 8.1] for a proof)

\medskip

\noindent (b) (Bernstein) The center $\zhkq$ of $\hkq$ is a
polynomial algebra over $\bbC$ in $n$-variables, generated by
$S_{x_i},\ i=1,2,...,n$. The elements $S_x,\ x\in X^+$, form a
$\bbC$-basis of the center $Z(\hkq)$ of $\hkq$.

\medskip

Note that $\Tq$ is a free $\zhkq$-module of rank $|W_0|$ (see [S]).

\medskip

\noindent{\bf 1.4.} Let $G$ be a simply connected simple algebraic
group over $\bbC$ with root system $R$. Let $T$ be a maximal torus
of $G$. Then we can identify $W_0$ with $N_G(T)/T$ and identify $X$
with Hom$(T,\bbC^*)$.

%\section{Isomorphism relations}

\def\hkp{H_{k,p}}
\def\zhkp{Z(H_{k,p}}
\def\hq{\hkq}
\def\hqs{H_{q,s}}
\def\a{\alpha}
\def\cws{\bbC[W]_s}
\def\cw{\bbC[W]}
\def\calN{\mathcal N}
\def\nqs{\calN_{q,s}}

Let $\mathcal C$ be a semisimple class of $G$. Choose an element $s$
in $\mathcal C\cap T$. The map $\theta_x\to x(s),\ x\in X$ defines a
homomorphism $\phi'_{q,s}:\ \Tq\to \bbC$. It is known that any
algebra homomorphism from $Z(H_q)$ to $\bbC$ is the restriction of
$\phi'_{q,s}$ for some $s\in T$. We shall denote by $\phi_{q,s}$ the
restriction to $Z(H_q)$ of $\phi'_{q,s}$. Let $I_{q,s}$ be the
two-sided ideal of $\hkq$ generated by all $S_x-\phi_{q,s}(S_x),\
x\in X^+$. Let $\hqs$ be the quotient algebra $\hq/I_{q,s}$. Then
$\hqs$ is a $\bbC$-algebra of dimension $|W_0|^2$. Since for any
$t\in\mathcal C\cap T$, we have $\hqs=H_{q,t}$, the algebra $\hqs$
depends only on the semisimple class $\mathcal C$.

We shall say that the central character $\phi_{q,s}$ admits
one-dimensional representations if $\hqs$  has one-dimensional
representations.

%If $s,t$ are not conjugate in $G$, then $\hqs$ and $H_{q,t}$ are
%different algebras.

\medskip

\noindent{1.5.}  Let $\mathfrak g$ be the Lie algebra of $G$ and
$\calN$ be the nilpotent cone of $\mathfrak g$. For any semisimple
element $s$ in $G$, define $\mathcal N_{s,q}$ to be the subset of
$\calN$ given by $\mathcal N_{s,q}=\{N\ |\ N\in \calN,\ {\text
Ad}(s)N=qN\}$. For a nilpotent element $N$ in $\calN_{s,q}$, let
$\cb^s_N$ be the variety consisting of all Borel subalgebras of
$\mathfrak g$ which contain $N$ and are fixed by ${\text Ad}(s)$.
Let $C_G(s)$ (resp. $C_G(N)$) be the centralizer of $s$ (resp. $N$)
in $G$. Denote by $A(s,N)$ the component group $C_G(s)\cap
C_G(N)/(C_G(s)\cap C_G(N))^o$ of $C_G(s)\cap C_G(N)$. Then $A(s,N)$
acts on the total complex coefficient Borel-Moore homology group
$H_*(\cb^s_N)$. Let $A(s,N)^\vee$ be the set of irreducible
representations of $A(s,N)$ that appear in $H_*(\cb^s_N)$.

The group $G$ acts on the set $G_{ss}$ of semisimple elements of $G$
by conjugacy and act on the variety $\calN$ by adjoint action. So
$G$ acts on the set of the union of all $A(s,N)^\vee$. If
$\sum_{w\in W_0}q^{l(w)}\ne 0$, the Deligne-Langlands-Lusztig
classification says that the isomorphism classes of irreducible
representations of $H_q$ is in one-to-one correspondence to the
$G$-orbits in the set of all triples $(s,N,\rho)$, $s\in G_{ss},\
N\in\calN$, $\rho\in A(s,N)^\vee$. See [K, KL2,X2].

 Let $s$ be a semisimple element in $G$. The group
$C_G(s)$ acts on the variety $\calN_{s,q}$ through the adjoint
action of $G$. So $C_G(s)$ acts on the set of all $A(s,N)^\vee$, $
N\in\nqs.$ The Deligne-Langlands-Lusztig classification is
equivalent to the following assertion.

\medskip

\noindent(a)If $\sum_{w\in W_0}q^{l(w)}\ne 0$, then the isomorphism
classes of irreducible representations of $\hqs$ is in one-to-one
correspondence to the $C_G(s)$-orbits of the pairs $(N,\rho)\
(N\in\nqs,\rho\in A(s,N)^\vee)$.

Since the variety $\cb_N^s$ is not empty for any $N\in \nqs$ [KL2],
by (a) we get

\medskip

\noindent(b) The number of isomorphism classes of irreducible
representations of $\hqs$ is not less than the number of the
$C_G(s)$-orbits in $\nqs$.

\medskip

\noindent{\bf1.6.} For a root $\alpha$ in $R$, let $e_\alpha$ be a
non zero element in the root subspace $\mathfrak g_\alpha$. Let
$\a_1,...,\alpha_n$ be the simple roots of $R$. We sometimes use
$e_i$ for $e_{\a_i}$.

For a root $\a$ in $R$, let $U_\a$ be the corresponding one
parameter subgroup of $G$ and $u_\a:\bbC\to U_\a$ be an isomorphism
such that $tu_{\a}(c)t^{-1}=u_\a(\a(t)c)$ for any $c\in\bbC$ and
$t\in T$.

\medskip

\centerline{\bf 2. Group algebra $\bbC[W]$  of $W$}

\medskip

In this section we consider the group algebra $\bbC[W]$, which is
identified with $H_1$.

\medskip

\noindent{\bf Lemma 2.1.}  Let $s$ be a semisimple element in $T$.

\noindent(a) If  $H_{1,s}$ has a one-dimensional representation,
then $s$ is in the center of $G$, i.e. $\a(s)=1$ for simple roots
$\a$ in $R$.

\noindent(b) If $s$ is in the center of $G$, then $H_{1,s}$ has a
one-dimensional representation.

Proof. (a) The following argument is provided by the referee.
Suppose that $\pi : H_1\to\bbC$ is a one-dimensional representation
of $H_1$ (and of $H_{1,s}$). This $\pi$ determines a character of $X
$, hence an element $t$ of $T$ conjugate to $s$, by $\pi(x) = x(t)$
for $x\in X$. Since $\pi(wxw^{-1}) = \pi(x)$ for $x\in X$, $w\in
W_0$, we see that $w(t) = t$ for any $w\in W_0$. Hence $t$ is in the
center of $G$ and $s = t$, (a) is true.

(b) Let $\Omega$ be the subgroup of $W$ consisting of elements of
length 0. Then $\Omega$ is a finite abelian group and is isomorphic
to the center of $G$. We identify the character group $\Omega^*$ of
$\Omega$ with the center  of $G$. If $s$ is an element in the center
of $G$ (which is identified with $\Omega^*)$, then the
one-dimensional representation $\chi:H_1\to\bbC$ defined by
$\chi(w)=1$ for all $w\in W_a$ and $\chi(\omega)=s(\omega)$ for any
$\omega\in\Omega$ has in fact central character $\phi_{1,s}$. Hence,
$H_{1,s}$ has a one-dimensional representation.

The lemma is proved.

\medskip

%Recall that for $s,t\in T$, the algebras $H_{1,s}$ and $H_{1,t}$ are
%the same if and only if $s$ and $t$ are conjugate in $G$.

 \noindent{\bf Corollary 2.2.} The number of central characters
 $\phi_{1,s}\ (s\in T)$  which admit  one-dimensional
representations is equal to the cardinality of the center of $G$.

\medskip

\noindent{\bf 2.3.} Assume that  $s$ is in the center. Then $A(s,N)$
is just the component group $A(N)=C_G(N)/C_G(N)^\circ$ of $C_G(N)$
for any nilpotent element $N$ in $\mathfrak g$. Any Borel subalgebra
of $\mathfrak g$ is fixed by Ad$(s)$. So the variety of $\cb_N^s$ is
just the variety $\cb_N$ of all Borel subalgebra of $\mathfrak g$
containing $N$. It is known that an irreducible representation of
$A(N)$ appears in the total Borel-Moore homology of $\cb_N$ if and
only if it appears in the top homology of the $\cb_N$.  Combining
the Deligne-Langlands-Lusztig classification for $H_1$ (see 1.5 (a))
and Springer's correspondence for $W_0$ we get the following fact.

\medskip

\noindent(a) If $s$ is in the center of $G$, then the  isomorphism
classes of irreducible representations of $H_{1,s}$ is in one-to-one
correspondence to the isomorphism classes of irreducible complex
representations of $W_0$.

\medskip

The number $|\text{Irr}(W_0)|$ of isomorphism classes of irreducible
complex representations of $W_0$ is well known. For convenience, we
list them as follows:

Type $A_n$: the number of all partitions of $n+1$,

Type $B_n$ $(n\ge 2)$ and $C_n$ $(n\ge 3)$: the number of ordered
pairs $(\xi,\eta)$ of partitions $\xi,\eta$ with $|\xi|+|\eta|=n$,

Type $D_n$ $(n\ge 4)$: the number of unordered pairs $(\xi,\eta)$ of
partitions $\xi,\eta$ with $|\xi|+|\eta|=n$, and any pairs
$(\xi,\eta)$ with $\xi=\eta$ and $|\xi|+|\eta|=n$ is counted twice.

Type $E_6,\ E_7,\ E_8$: 25, 60, 112,   respectively

Type $F_4$: 25,

Type $G_2$: 6.

\medskip

\centerline{\bf 3. Some particular semisimple elements in $T$}

\medskip

In this section we assume that $q\ne 1$. We are interested in the
quotient algebras $\hqs$ which have one-dimensional representations.

\medskip

\noindent {\bf Lemma 3.1.}  Assume that $\hqs$ has a one-dimensional
representation  on which all $T_r$ ($r\in S$) act by the same scalar
multiplication. Then in the conjugacy class of $s$ we can find an
element $t\in T$ such that $\alpha(t)=q$ for all simple roots
$\alpha$.

% Assume that $s\in T$ and $\a(s)=q$ for all simple roots $\a$.
%Then $\hqs$ has a one-dimensional representation on which all $T_r$
%($r\in S$) act by the same scalar multiplication.

Proof.  The scalar is $q$ or $-1$.  Using the explicit formulas for
reduced expressions of fundamental weights in [L2] we can see that
the lemma is true. We illustrate the argument by using type $\tilde
A_n$.

 Assume that $R$ is of type $A_n$. We number the simple reflections $r_0,r_1,...,r_n$
 as usual and simply write $T_i$ for $T_{r_i}$.
 Let $\a_1,...,\a_n$ be the simple roots and $x_1,...,x_n$ be the
 corresponding fundamental weights of the weight lattice $X$. According to [L2],
 for $1\le i \le n$, we have
$$T_{x_i}=T_{\tau^{n+1-i}}(T_{n+1-i}T_{n+2-i}\cdots
T_n)(T_{n-i}T_{n+1-i}\cdots T_{n-1})
 \cdots (T_1T_2\cdots T_i),$$
 where $\tau\in W$ has length 0 and $\tau r_0=r_{1}\tau$, $\tau r_1=r_{2}\tau$,..., $\tau
 r_n=r_{0}\tau$. Note that $\tau^{n+1}=e$ is the neutral element of
 $W$.

 In $X$ we set $x_0=x_{n+1}=0$. Then the length of $x_i$ is $i(n+1-i)$. By
 definition, for $1\le i\le n$, we have
 $$\theta_{\a_i}=(q^{-i(n+1-i)}T_{x_i}^2)(q^{-(i+1)(n-i)/2}T_{x_{i+1}})^{-1}(q^{-(i-1)(n+2-i)/2}T_{x_{i-1}})^{-1}.$$
 If all $T_i$ act on a one-dimensional representation of $\hqs$ by
 scalar $q$, then for all $1\le i\le n$,  $\theta_{\a_i}$ act on
 the one-dimensional representation of $\hqs$ by
 scalar $q$. So the required $t$ exists in this case.

 If all $T_i$ act on a one-dimensional representation of $\hqs$ by
 scalar $-1$, then for all $1\le i\le n$,  $\theta_{\a_i}$ act on
 the one-dimensional representation of $\hqs$ by
 scalar $q^{-1}$. This means that in the conjugacy class $\mathcal C$ of $s$, there exists $t'\in \mathcal C\cap T$ such that
 $\a_i(t')=q^{-1}$ for all $i=1,2,...,n$. Let $w_0$ be the longest element of $W_0$ and Set $t=w_0(t')$. Then
  $t$ is in $\mathcal C\cap T$ and $\a_i(t)=q$ for all
 $i=1,2,...,n$. The lemma is proved for type $\tilde A$. The proof for other types is similar.

\medskip

\def\ve{\varepsilon}

\noindent{\bf Lemma 3.2.} Assume that $(q-1)\sum_{w\in
W_0}q^{l(w)}\ne 0$ and  $\hqs$ has a one-dimensional representation
on which  some $T_r$ ($r\in S$) act by  scalar multiplication of $q$
and some $T_r$ ($r\in S$) act by  scalar multiplication of $-1$.
Then in the conjugacy class of $s$ there is no  element $t\in T$
such that $\alpha(t)=q$ for all simple roots $\alpha$.

Proof.  The root system $R$ must be one of the following types:
$B_n\ (n\ge 2)$, $C_n\ (n\ge 3)$, $F_4$, $G_2$. We prove the lemma
case by case. Let $t\in T$ be such that $\a(t)=q$ for all simple
roots $\a$.

\noindent{\bf Type $B_n\ (n\ge 2)$.}  There exist
$\varepsilon_1,\ve_2,...,\ve_n$ in Hom$(T,\bbC^*)$ such that
$\a_i=\ve_i-\ve_{i+1}$ for $i=1,...,n-1$ and $\a_n=\ve_n$. The
maximal exponent of $R$ is $2n-1$ and $\sum_{w\in
W_0}q^{l(w)}=(q-1)^{-n}\prod_{i=1}^n(q^{2i}-1)$.

Let $r_i$ be the simple reflection corresponding to $\a_i$ and $r_0$
the simple reflection out of $W_0$. Assume
$T_{r_1},T_{r_2},...,T_{r_{n-1}}$ act on the one-dimensional
representation by scalar $q$ and $T_{r_0}$, $T_{r_n}$ act on it by
scalar $-1$. Using the explicit formula for $\theta_{x_i}$ in [L2]
we see that $\a_i(s)=q$ for $i=1,2,...,n-1$ and $\a_n(s)=q^{-1}$.
Therefore $C_G(s)$ contains $T, U_{\pm\ve_{n-1}},
U_{\pm(\ve_{n-2}+\ve_n)}$.

 If
the order $o(q)$ is greater than the maximal exponent $2n-1$, then
$t$ is regular and $C_G(t)=T$. In this case $s$ and $t$ are not
conjugate in $G$.

Now assume that $o(q)\le 2n-1$. Our assumption on $q$ implies that
$o(q)$ is odd and $n+1\le o(q)\le 2n-1$. Let $o(q)=n+i$ for some
$1\le i\le n-1$. Since $o(q)=n+i$, The centralizer $C_G(t)$ of $t$
is generated by $T$ and all $U_\a$ with $\a(t)=1$. It is easy to see
that $\a(t)=1$ if and only if $\a$ is one of the following roots:
$\pm(\ve_j+\ve_{n+2-j-i})$, $1\le j< \frac{n+2-i}2$. Noting that all
the roots $\ve_j+\ve_{n+1-j-i}$ are long roots and $\ve_{n-1}$ is
short root, we see that $C_G(s)$ and $C_G(t)$ are not conjugate in
$G$, so $s$ and $t$ are not conjugate in $G$.

 Assume
$T_{r_1},T_{r_2},...,T_{r_{n-1}}$ act on the one-dimensional
representation by scalar $-1$ and $T_{r_0}$, $T_{r_n}$ act on it by
scalar $q$. Again using the explicit formula for $\theta_{x_i}$ in
[L2] we see that $\a_i(s)=q^{-1}$ for $i=1,2,...,n-1$ and
$\a_n(s)=q$. Let $g\in N_G(T)$ be a representative of the longest
element of $W_0$ and $s'=gsg^{-1}$. Then $\a_i(s')=q$ for
$i=1,2,...,n-1$ and $\a_n(s')=q^{-1}$. By the above argument we see
that $s'$ and $t$ are not conjugate in $G$.

\noindent{\bf Type $C_n$ ($n\ge 3)$.} The argument is similar to
that for  type $B_n$.

\noindent{\bf Type $F_4$.} Let $V=\mathbb{R}^4$ and
$\ve_1,...,\ve_4$ the standard basis of $V$. We may assume that $R$
consists of the following elements:

$\pm\ve_i\ (1\le i\le 4)$, \quad $\pm\ve_i\pm\ve_j\ (1\le i<j\le
4)$,

$\frac12(\pm\ve_1\pm\ve_2\pm\ve_3\pm\ve_4)$.

The maximal exponent of $R$ is 11 and $$\sum_{w\in
W_0}q^{l(w)}=(q-1)^{-3}(q+1)(q^6-1)(q^8-1)(q^{12}-1).$$

Let $r_i$ be the simple reflection corresponding to $\a_i$ and $r_0$
is the reflection out of $W_0$. Then $r_0r_4r_0=r_4r_0r_4$. Assume
$T_{r_1},T_{r_2}$ act on the one-dimensional representation by
scalar $q$ and $T_{r_3},T_{r_4}, T_{r_0}$ act on it by scalar $-1$.

According to [L2, p.646], in $H_q$ we have
$$\theta_{x_4}=q^{-8}T_{r_0}T_{r_4}T_{r_3}T_{r_2}T_{r_1}T_{r_3}T_{r_4}T_{r_2}T_{r_3}T_{r_2}T_{r_4}T_{r_3}
T_{r_1}T_{r_2}T_{r_3}T_{r_4},$$
$$\theta_{x_3}\theta_{x_4}^{-1}=qT_{r_4}^{-1}\theta_{x_4}T_{r_4}^{-1},$$
$$\theta_{x_2}\theta_{x_3}^{-1}=qT_{r_3}^{-1}\theta_{x_3}\theta_{x_4}^{-1}T_{r_3}^{-1},$$
$$\theta_{x_1}\theta_{x_2}^{-1}\theta_{x_3}=qT_{r_2}^{-1}\theta_{x_2}\theta_{x_3}^{-1}T_{r_2}^{-1}.$$

So we have $\a_1(s)=\a_2(s)=q$ and $\a_3(s)=\a_4(s)=q^{-1}$.  The
centralizer $C_G(s)$ of $s$ contains $T$, $U_{\pm(\a_2+\a_3)},$
$U_{\pm(\a_1+\a_2+2\a_3)},$ $U_{\pm(\a_1+\a_2+\a_3+\a_4)},$
$U_{\pm(\a_1+2\a_2+2\a_3+\a_4)}.$

If the order $o(q)$ is greater than 11, then $t$ is regular and
$C_G(s)=T$, so $s$ and $t$ are not conjugate.

Now assume that $o(q)\le 11$. Our assumption on $q$ implies that
$o(q)=5,7,9,10,11$. If $o(q)=11$, then $C_G(t)$ is generated by $T$
and $U_{\pm\beta}$, where $\beta$ is the highest root in $R$. If
$o(q)=10$, then $C_G(t)$ is generated by $T$ and
$U_{\pm(\beta-\a_1)}$. If $o(q)=9$, then $C_G(t)$ is generated by
$T$, $U_{\pm(\beta-\a_1-\a_2)}$. If $o(q)=7$, then
 $C_G(t)$ is generated by
$T$, $U_{\pm(\beta-\a_1-\a_2-2\a_3)}$.
$U_{\pm(\beta-\a_1-\a_2-\a_3-\a_4)}$. If $o(q)=5$, then $C_G(t)$ is
generated by $T$, $U_{\pm(\a_1+2\a_2+2\a_3)}$,
$U_{\pm(\a_2+2\a_3+2\a_4)}$, $U_{\pm(\beta-\a_1)}$. In all cases
$C_G(s)$ and $C_G(t)$ are not conjugate in $G$. Hence $s$ and $t$
are not conjugate in $G$.

Assume $T_{r_1},T_{r_2}$ act on the one-dimensional representation
by scalar $-1$ and $T_{r_3},T_{r_4}, T_{r_0}$ act on it by scalar
$q$. One can check as above that $\a_1(s)=\a_2(s)=q^{-1}$ and
$\a_3(s)=\a_4(s)=q$. Let $g\in N_G(T)$ be a representative of the
longest element of $W_0$ and $s'=gsg^{-1}$. Then
$\a_1(s')=\a_2(s')=q$ and $\a_3(s')=\a_4(s')=q^{-1}$. By the above
discussion we know that $s'$ and $t$ are not conjugate in $G$.

\noindent{\bf Type $G_2$.}  We number the simple reflections
$r_0,r_2,r_2$ so that $r_0r_2=r_2r_0$. In $W_0$ we have
$x_1=r_0r_1r_2r_1r_2r_1$ and $x_2=r_0r_1r_2r_1r_2r_0r_1r_2r_1r_2$.
In the weight lattice we have $x_1=2\a_1+\a_2$ and
$x_2=3\a_1+2\a_2$, where $\a_i$ are the simple roots.  If $T_{r_1}$
and $T_{r_0}$ act on the one-dimensional representation by scalar
$q$ and $T_{r_2}$ acts it by scalar $-1$, then both $\theta_{x_1}$
and $\theta_{x_2}$ act on it by scalar $q$. So $\a_1(s)=q$ and
$\a_2(s)=q^{-1}$. Note that $(\a_1+\a_2)(s)=1$. So $C_G(s)$ contains
$T$ and $U_{\pm(\a_1+\a_2)}$.

 If the order $o(q)$
of $q$ greater than 5, then $C_G(s)=T$. So $s$ and $t$ are not
conjugate. If $o(q)\le 5$, then $o(q)=5,4$ since $\sum_{w\in
W_0}q^{l(w)}\ne 0$. In these cases $C_G(s)$ is generated by $T$ and
$U_{\pm(\a_1+\a_2)}$, $C_G(t)$ is generated by $T$ and
$U_{\pm(3\a_1+\a_2)}$ if $o(q)=4$ or $U_{\pm(3\a_1+2\a_2)}$ if
$o(q)=5$. Clearly $C_G(s)$ and $C_G(t)$ are not conjugate, so $s$
and $t$ are not conjugate.

If $T_{r_1}$ and $T_{r_0}$ act on the one-dimensional representation
by scalar $-1$ and $T_{r_2}$ acts it by scalar $q$, then
$\a_1(s)=q^{-1}$ and $\a_2(s)=q$. Let $g\in N_G(T)$ be a
representative of the longest element of $W_0$ and $s'=gsg^{-1}$.
Then $\a_1(s')=q$ and $\a_2(s')=q^{-1}$. By the above discussion we
know that $s'$ and $t$ are not conjugate in $G$.

 The lemma is proved.

\medskip

\noindent{\bf Corollary 3.3.} Assume that the root system $R$ is not
simply laced. If $(q-1)\sum_{w\in W_0}q^{l(w)}\ne 0$, then the
 number of central characters
 $\phi_{q,s}\ (s\in T)$  which admit  one-dimensional
representations is twice the cardinality of the center of $G$.

Proof.  Let $s\in T$ be such that $\a(s)=q$ for all simple roots
$\a$. It is easy to check that the map $T_r\theta_x\to qx(s)$ for
any simple reflection $r$ in $W_0$ and $x\in X$ defines an algebra
homomorphism $H_q\to \mathbb{C}$. Hence, the central character
$\phi_{q,s}$ admits one-dimensional representations of $H_q$ on
which all $T_i$ act by the same scalar. Since $(q-1)\sum_{w\in
W_0}q^{l(w)}\ne 0$, we have $q\ne -1$. So we can find $t\in T$ such
that the central character $\phi_{q,t}$ admits one-dimensional
representations of $H_q$ on which some $T_i$ act by
 scalar $q$ and some $T_i$ act by scalar $-1$. By the calculations in the proof of Lemma 3.2, it is no harm to assume that
 $\a(t)=q$ for short simple roots $\a$ and $\a(t)=q^{-1}$ for long  simple roots $\a$.

By Lemma 3.2, $s$ and $t$ are not conjugate in $G$ if
$(q-1)\sum_{w\in W_0}q^{l(w)}\ne 0$.

 For types $F_4$ and
$G_2$, the center of $G$ is trivial and the weight lattice equals
the root lattice. Hence there are exactly two central characters of
$H_q$ which admit  one-dimensional representations if
$(q-1)\sum_{w\in W_0}q^{l(w)}\ne 0$.

Now assume that $R$ is of type $B_n$ ($n\ge 2)$ or $C_n$ ($n\ge 3)$.
Let $s'\in T$. Assume that the central character $\phi_{q,s'}$
admits one-dimensional representations of $H_q$ on which all $T_i$
act by the same scalar. By Lemma 3.1, in the conjugacy class of $s'$
we can find $s''\in T$ such that $\a(s'')=q$ for all simple roots
$\a$. Then $s''s^{-1}=c$ is in the center of $G$. We need show that
$s''$ and $s$ are not conjugate in $G$ if $c$ is non-trivial and
$(q-1)\sum_{w\in W_0}q^{l(w)}\ne 0$.

First consider type $B_n$. In this case $G=Spin_{2n+1}(\bbC)$ and
$SO_{2n+1}(\bbC)$ is its quotient group. Clearly,  $s''$ and $s$
have the same image in $SO_{2n+1}(\bbC)$. It is no harm to assume
the image is
$$\text{diag}(1,q^{n},...,q^{2},q^{1},q^{-n},...,q^{-2},q^{-1}).$$
If $s''$ and $s$ are conjugate in $G$, we must have $q=q^i$ for some
$2\le i\le n$ or $q=q^{-j}$ for some $1\le j\le n$. If $j\ne n$,
then $(q-1)\sum_{w\in W_0}q^{l(w)}=0$. If $j=n$, then $q^{n+1}=1$.
We also have $q^{\frac{n(n+1)}2}=-1$ in this case since $c$ is
non-trivial. This forces that $n+1$ is even and $\sum_{w\in
W_0}q^{l(w)}=0$. So $s''$ and $s$ are not conjugate if
$(q-1)\sum_{w\in W_0}q^{l(w)}\ne 0$.

Now assume that $R$ is of type $C_n$. Then $G=Sp_{2n}(\bbC)$. If $c$
is non-trivial, then we may assume that $s,s''$ are the following
two elements:
$$\text{diag}(q^{\frac{2n-1}2},...,q^{\frac32},q^{\frac12},q^{-\frac{2n-1}2},...,q^{-\frac32},q^{-\frac12})$$
$$\text{diag}(-q^{\frac{2n-1}2},...,-q^{\frac32},-q^{\frac12},-q^{-\frac{2n-1}2},...,-q^{-\frac32},-q^{-\frac12}).$$
If $s''$ and $s$ are conjugate, then $q^{\frac12}=-q^{\frac12+i}$
for some $n-1\ge i\ge 1$ or $q^{\frac12}=-q^{-\frac12-j}$ for some
$n-1\ge j\ge 0$. So, $q^i=-1$ or $q^{1+j}=-1$. This contradicts that
$(q-1)\sum_{w\in W_0}q^{l(w)}\ne 0$. Hence,  $s''$ and $s$ are not
conjugate in $G$.

Assume that the central character $\phi_{q,t'}$ admits
one-dimensional representations of $H_q$ on which some $T_i$ act by
the scalar $q$ and some $T_i$ act by the scalar $-1$. By the
calculations in the proof of Lemma 3.2, in the conjugacy class of
$t'$ we can find $t''\in T$ such that $\a(t'')=q$ for all simple
short roots $\a$ and $\a(t'')=q^{-1}$ for all simple long roots
$\a$. Then $t''t^{-1}=c$ is in the center of $G$. As above, we can
show  that $t''$ and $t$ are not conjugate in $G$ if $c$ is
non-trivial and $(q-1)\sum_{w\in W_0}q^{l(w)}\ne 0$.

The corollary is proved.

\medskip

\noindent{\bf Lemma 3.4.} Assume that the order of $q$ is greater
than the maximal exponent of the root system $R$. Let $s\in T$ be
such that $\a(s)=q$ for all simple roots $\a$ in $R$. Then

\noindent(a) The number of $C_G(s)$-orbits of the variety
$\calN_{q,s}$ is equal to $2^n$.

\noindent(b) For any $N\in \nqs$, the group $A(s,N)$ acts on
$H_*(\cb^s_N)$ trivially, so $A(s,N)^\vee$ contains only the trivial
representation of $A(s,N)$.

\noindent(c) The number of the isomorphism classes of irreducible
representations of $\hqs$ is $2^n$.

Proof.  Since the order of $q$ is greater than the maximal exponent
of the root system $R$, we see that $s$ is regular, so $C_G(s)=T$.
Moreover, we have $\calN_{q,s}=\{\sum_{1\le i\le n}a_ie_i\,|\,
a_i\in\bbC\}.$ For a  subset $I$ of $\{1,2,...,n\}$, set
$e_I=\sum_{i\in I}e_{i}.$ (We understand that $e_I=0$ if $I$ is the
empty set.) Then any element in $\calN_{q,s}$ is $C_G(s)$-conjugate
to some $e_I$. Clearly if $I\ne J$, the $e_I$ and $e_J$ are not in
the same $C_G(s)$-orbit. Part (a) is proved.

Since $s$ is regular, the variety $\cb^s$ of Borel subalgebras of
$\mathfrak g$ fixed by Ad$(s)$ is exactly the variety $\cb^T$ of
Borel subalgebras of $\mathfrak g$ fixed by Ad$T$. Now $C_G(s)=T$,
so for any $N\in\nqs$, the group $C_G(s)\cap C_G(N)$ acts trivially
on $\cb^s_N$. Part (b) follows.

Part (c) follows from (a) and (b). The lemma is proved.

\medskip

\noindent{\bf 3.5.} In the rest part of this section $R$ is assumed
to be simply laced and $\hqs$ is required to satisfy the assumption
of Lemma 3.1. We hope to know the number of the $C_G(s)$-orbits in
$\nqs$, provided that $(q-1)\sum_{w\in W_0}q^{l(w)}\ne 0$ and the
order $o(q)$ of $q$ is not greater than the maximal exponent of $R$.
The following fact of Lusztig [L3] will be useful.

\medskip

\noindent(a) Assume  $x\in \mathfrak g$ and Ad$(s)(x)=qx$. If
$(q-1)\sum_{w\in W_0}q^{l(w)}\ne 0$, then $x\in\nqs$.

\medskip

 Since $\a(s)=q$ for all simple roots $\a$ in $R$,
by (a) we get

\medskip

\noindent(b) $\nqs$ is a linear space spanned by all $e_\beta$ with
$\beta(s)=q$ and is naturally a $C_G(s)$-module.

Since $T\subset C_G(s)$, any $C_G(s)$-submodule of $\nqs$ is spanned
by some $e_\beta$ with $\beta(s)=q$. Since $C_G(s)$ is reductive,
$\nqs$ is the direct sum of some irreducible $C_G(s)$-submodules of
$\nqs$.

\medskip

We discuss, on a case by case basis,  the number of $C_G(s)$-orbits
in $\nqs$. To do this, we first decompose $\nqs$ into the direct sum
of some irreducible $C_G(s)$-submodules. Let $\beta\in R$ be such
that $\beta(s)=1$ and $u_\beta(\xi)$ ($\xi\in \bbC$, see 1.6 for
definition) an element in the one parameter subgroup $U_\beta$ of
$G$. Then we may assume that $u_\beta(\xi)$ acts on $\mathfrak g$
(hence $\nqs$) through $\exp(\xi\text{ad}(e_\beta))$. This idea is
useful in determining the submodule structure of $\nqs$. Recall that
we identify the weight lattice $X$ of $R$ with Hom$(T,\bbC^*)$.

\medskip

\noindent{\bf Type $A_n$} The maximal exponent is $n$ and
$\sum_{w\in W_0}q^{l(w)}=(q-1)^{-n}\prod_{i=2}^{n+1}(q^i-1)$. There
are no complex numbers $q\ne 1$ such that its order is not greater
than $n$ and $\sum_{w\in W_0}q^{l(w)}\ne 0$.

\medskip

\noindent{\bf Type $D_n$ ($n\ge 4)$}  There exist
$\varepsilon_1,\ve_2,...,\ve_n$ in Hom$(T,\bbC^*)$ such that
$\a_i=\ve_i-\ve_{i+1}$ for $i=1,...,n-1$ and $\a_n=\ve_{n-1}+\ve_n$.
The maximal exponent of $R$ is $2n-3$ and  $\sum_{w\in
W_0}q^{l(w)}=(q-1)^{-n}(q^n-1)\prod_{i=1}^{n-1}(q^{2i}-1)$. The
assumption on $q$ implies that $o(q)$ is odd and $n+1\le o(q)\le
2n-3$. Let $o(q)=n+i$ for some $1\le i\le n-3$.

By assumption, $\a_i(s)=q$. Since $o(q)=n+i$, we see the other roots
satisfying $\a(s)=q$ are the following: $\ve_j+\ve_{n-1-j-i}$,
$-\ve_k-\ve_{n+1-k-i}$,   with $1\le j< \frac{n-i-1}2$ and $1\le k<
\frac{n+1-i}2$.

It is easy to check that $\a(s)=1$ if and only if $\a$ is one of the
following roots: $\pm(\ve_j+\ve_{n-j-i})$, $1\le j< \frac{n-i}2$.

If $i=n-3$, then $C_G(s)$ is generated by $T$ and
$U_{\pm(\ve_1+\ve_2)}$. The linear space $\nqs$ is spanned by
$e_{\ve_1-\ve_2},$ $e_{\ve_2-\ve_3},$,..., $e_{\ve_{n-1}-\ve_{n}}$,
$e_{\ve_{n-1}+\ve_n}$ and $e_{-\ve_1-\ve_3}$. It is easy to see that
$e_{\ve_2-\ve_3}$ and $e_{-\ve_1-\ve_3}$ span a $C_G(s)$-submodule
$M_2$ of $\nqs$. For any $1\le j\ne 2\le n-1$, the element
$e_{\ve_{j}-\ve_{j+1}}$ spans a one-dimensional $C_G(s)$-submodule
$M_j$ of $\nqs$. The space $\nqs$ is the direct sum of the
submodules $M_1,M_2,...,M_n$ and each $M_j$ has two $C_G(s)$-orbits.
Therfore the number of $C_G(s)$-orbits in $\nqs$ is $2^n$ .

Now assume that $i\le n-5$. The following elements span a
$C_G(s)$-submodule $M_j$ of $\calN_{q,s}$: $e_{\ve_j-\ve_{j+1}},$ $
e_{-\ve_{j+1}-\ve_{n-j-i}}$, $e_{\ve_j+\ve_{n-j-1-i}}$,
$e_{\ve_{n-j-1-i}-\ve_{n-j-i}}$ for $j=1,2,...,\frac {n-i-3}2$. The
elements $e_{\ve_{j}-\ve_{j+1}}$ and $e_{-\ve_{j-1}-\ve_{j+1}}$ span
a two-dimensional $C_G(s)$-submodule $M'_j$ of $\calN_{q,s}$ for
$j=\frac{n-i+1}2$. The element $e_{\ve_{j}-\ve_{j+1}}$ spans a
one-dimensional $C_G(s)$-submodule $M'_{j}$ of $\calN_{q,s}$ for
$j=\frac{n-i-1}2, n-i,n-i+1,...,n-1$. Also the element
$e_{\ve_{n-1}+\ve_n}$  spans  a one-dimensional $C_G(s)$-submodule
$M'_{n}$ of $\calN_{q,s}$ Clearly $\calN_{q,s}$ is the direct sum of
all the submodules $M_j, M'_k$. One may check easily that each $M_j$
has three $C_G(s)$-orbits and each $M'_k$ has two $C_G(s)$-orbits.
So the number of $C_G(s)$-orbits in $\calN_{q,s}$ is not less than
$2^{i+3}\cdot 3^{\frac {n-i-3}2}$.

\medskip

\noindent{\bf Type $E_6$} We number the simple roots as in [B]. One
has $$\sum_{w\in
W_0}q^{l(w)}=\frac{(q^2-1)(q^5-1)(q^6-1)(q^8-1)(q^9-1)(q^{12}-1)}{(q-1)^{6}}.$$
The maximal exponent is 11. The assumption on $q$ implies that
$o(q)=7,10,11$.

If $o(q)=7$, then $C_G(s)$ is generated by $T, U_{\pm\beta_i}$,
$i=1,2,3$, here $\beta_1=\a_1+\a_2+2\a_3+2\a_4+\a_5$,
$\beta_2=\a_2+\a_3+2\a_4+2\a_5+\a_6$,
$\beta_3=\a_1+\a_2+\a_3+2\a_4+\a_5+\a_6.$  If $\a$ is not simple
root, then $\a(s)=q$ if and only if $\a$ is one of the following
roots: $\gamma_1=\a_1+\a_2+2\a_3+2\a_4+\a_5+\a_6$,
$\gamma_2=\a_1+\a_2+\a_3+2\a_4+2\a_5+\a_6$,
$-\gamma_3=-(\a_1+\a_2+\a_3+2\a_4+\a_5)$,
$-\gamma_4=-(\a_2+\a_3+2\a_4+\a_5+\a_6)$,
$-\gamma_5=-(\a_1+\a_2+\a_3+\a_4+\a_5+\a_6)$.  It is easy to check
the following facts:

(1) $e_{\a_1}$, $e_{\gamma_2}, $ $e_{\a_5}$, $e_{-\gamma_4}$ span a
$C_G(s)$-submodule $M_1$ of $\nqs$,

(2) $e_{\a_3}$, $e_{\gamma_1}, $ $e_{\a_6}$, $e_{-\gamma_3}$ span a
$C_G(s)$-submodule $M_2$ of $\nqs$,

(3) $e_{\a_4},e_{-\gamma_5}$ span a $C_G(s)$-submodule $M_3$ of
$\nqs$,

(4) $e_{\a_2}$ spans a $C_G(s)$-submodule $M_4$ of $\nqs$.

Clearly $\nqs$ is the direct sum of all $M_i$. Noting that $M_i$ has
three $C_G(s)$-orbits for $i=1,2$ and $M_j$ has two $C_G(s)$-orbits
for $j=3,4$, we see that the number of $C_G(s)$-orbits in
$\calN_{q,s}$ is not less than 36.

If $o(q)=10$ or 11, then we can see easily that the number of
$C_G(s)$-orbits in $\calN_{q,s}$ is $2^6$.

\medskip

\noindent{\bf Type $E_7$} We number the simple roots as in [B]. One
has $$\sum_{w\in
W_0}q^{l(w)}=\frac{(q^2-1)(q^6-1)(q^8-1)(q^{10}-1)(q^{12}-1)(q^{14}-1)(q^{18}-1)}{(q-1)^{7}}.$$
The maximal exponent is 17. The assumption on $q$ implies that
$o(q)=11,13,15,16,17$.

If $o(q)=11$, then $C_G(s)$ is generated by $T, U_{\pm\beta_i}$,
$i=1,2,3$, here $\beta_1=\a_1+2\a_2+2\a_3+3\a_4+2\a_5+\a_6$,
$\beta_2=\a_1+\a_2+2\a_3+2\a_4+2\a_5+2\a_6+\a_7$,
$\beta_3=\a_1+\a_2+2\a_3+3\a_4+2\a_5+\a_6+\a_7.$ If $\a$ is not
simple root, then $\a(s)=q$ if and only if $\a$ is one of the
following roots: $\gamma_1=\a_1+\a_2+2\a_3+3\a_4+2\a_5+2\a_6+\a_7$,
$\gamma_2=\a_1+2\a_2+2\a_3+3\a_4+2\a_5+\a_6+\a_7$,
$-\gamma_3=-(\a_1+\a_2+2\a_3+3\a_4+2\a_5+\a_6)$,
$-\gamma_4=-(\a_1+\a_2+2\a_3+2\a_4+2\a_5+\a_6+\a_7)$,
$-\gamma_5=-(\a_1+\a_2+\a_3+2\a_4+2\a_5+2\a_6+\a_7)$.  It is easy to
check the following facts:

(5) $e_{\a_2}$, $e_{\gamma_2}, $ $e_{\a_7}$, $e_{-\gamma_4}$ span a
$C_G(s)$-submodule $M_2$ of $\nqs$,

(6) $e_{\a_4}$, $e_{\gamma_1}, $ $e_{\a_6}$, $e_{-\gamma_4}$ span a
$C_G(s)$-submodule $M_4$ of $\nqs$,

(7) $e_{\a_3},e_{-\gamma_5}$ span a $C_G(s)$-submodule $M_3$ of
$\nqs$,

(8) For $i=1$ or $5$, the element $e_{\a_i}$ spans a
$C_G(s)$-submodule $M_i$ of $\nqs$.

Clearly $\nqs$ is the direct sum of all $M_i$. Noting that $M_i$ has
three $C_G(s)$-orbits for $i=2,4$ and $M_j$ has two $C_G(s)$-orbits
for $j=1,3,5$, we see that the number of $C_G(s)$-orbits in
$\calN_{q,s}$ is not less than $2^3\times 3^2=72$.

If $o(q)=13$,  then $C_G(s)$ is generated by $T, U_{\pm\sigma_i}$,
$i=1,2$, here $\sigma_1=\a_1+\a_2+2\a_3+3\a_4+3\a_5+2\a_6+\a_7$,
$\sigma_2=\a_1+2\a_2+2\a_3+3\a_4+2\a_5+2\a_6+\a_7$. If $\a$ is not
simple root, then $\a(s)=q$ if and only if $\a$ is one of the
following roots: $\tau_1=\a_1+2\a_2+2\a_3+3\a_4+3\a_5+2\a_6+\a_7$,
$-\tau_2=-(\a_1+\a_2+2\a_3+3\a_4+2\a_5+2\a_6+\a_7)$,
$-\tau_3=-(\a_1+2\a_2+2\a_3+3\a_4+2\a_5+\a_6+\a_7)$.  It is easy to
check the following facts:

(9) $e_{\a_2}$, $e_{\tau_1}, $ $e_{\a_5}$, $e_{-\tau_2}$ span a
$C_G(s)$-submodule $M_2$ of $\nqs$,

(10) $e_{\a_6},e_{-\tau_3}$ span a $C_G(s)$-submodule $M_6$ of
$\nqs$,

(11) The element $e_{\a_i}$ spans a $C_G(s)$-submodule $M_i$ of
$\nqs$ for $i=1,3,4,7$.

Clearly $\nqs$ is the direct sum of all $M_i$. Noting that $M_2$ has
three $C_G(s)$-orbits and $M_j$ has two $C_G(s)$-orbits for
$j=1,3,4,6,7$, we see that the number of $C_G(s)$-orbits in
$\calN_{q,s}$ is not less than $2^5\times 3=96$.

If $o(q)=15,16,17$, then we can see easily that the number of
$C_G(s)$-orbits in $\calN_{q,s}$ is $2^7$.

\medskip

\noindent{\bf Type $E_8$} This is the most complicated case. We
number the simple roots as in [B]. One has $$\sum_{w\in
W_0}q^{l(w)}=\frac{(q^2-1)(q^8-1)(q^{12}-1)(q^{14}-1)(q^{18}-1)(q^{20}-1)(q^{24}-1)(q^{30}-1)}{(q-1)^{8}}.$$
The maximal exponent is 29. The assumption on $q$ implies that
$o(q)=11,13,16,17,19,21,22,23,25,26,27,28,29$.

If $o(q)=11$, then $C_G(s)$ is generated by $T, U_{\pm\beta_i}$,
$1\le i\le 6$, here

$\beta_1=\a_1+2\a_2+2\a_3+3\a_4+2\a_5+\a_6$,

$\beta_2=\a_1+\a_2+2\a_3+2\a_4+2\a_5+2\a_6+\a_7$,

$\beta_3=\a_1+\a_2+2\a_3+3\a_4+2\a_5+\a_6+\a_7,$

$\beta_4=\a_1+\a_2+2\a_3+2\a_4+2\a_5+\a_6+\a_7+\a_8,$

$\beta_5=\a_2+\a_3+2\a_4+2\a_5+2\a_6+2\a_7+\a_8,$

$\beta_6=\a_1+\a_2+\a_3+2\a_4+2\a_5+2\a_6+\a_7+\a_8,$

$\beta_7=\beta_1+\beta_5=\a_1+3\a_2+3\a_3+5\a_4+4\a_5+3\a_6+2\a_7+\a_8,$

$\beta_8=\beta_3+\beta_6=2\a_1+2\a_2+3\a_3+5\a_4+4\a_5+3\a_6+2\a_7+\a_8.$

 If $\a$ is not
simple root, then $\a(s)=q$ if and only if $\a$ is one of the
following roots:

$\gamma_1=\a_1+\a_2+2\a_3+3\a_4+2\a_5+2\a_6+\a_7$,

$\gamma_2=\a_1+2\a_2+2\a_3+3\a_4+2\a_5+\a_6+\a_7$,

$-\gamma_3=-(\a_1+\a_2+2\a_3+3\a_4+2\a_5+\a_6)$,

$-\gamma_4=-(\a_1+\a_2+2\a_3+2\a_4+2\a_5+\a_6+\a_7)$,

$-\gamma_5=-(\a_1+\a_2+\a_3+2\a_4+2\a_5+2\a_6+\a_7)$,

$-\gamma_6=-(\a_2+\a_3+2\a_4+2\a_5+2\a_6+\a_7+\a_8)$,

$-\gamma_7=-(\a_1+\a_2+2\a_3+2\a_4+\a_5+2\a_6+\a_7+\a_8)$,

$-\gamma_8=-(\a_1+\a_2+\a_3+2\a_4+2\a_5+\a_6+\a_7+\a_8)$,

$\gamma_9=\a_1+\a_2+2\a_3+3\a_4+2\a_5+\a_6+\a_7+\a_8$,

$\gamma_{10}=\a_1+\a_2+2\a_3+2\a_4+2\a_5+2\a_6+\a_7+\a_8$,

$\gamma_{11}=\a_1+\a_2+\a_3+2\a_4+2\a_5+2\a_6+2\a_7+\a_8$,

$\gamma_{12}=2\a_1+3\a_2+3\a_3+5\a_4+4\a_5+3\a_6+2\a_7+\a_8$,

$\gamma_{13}=2\a_1+2\a_2+4\a_3+5\a_4+4\a_5+3\a_6+2\a_7+\a_8$,

$-\gamma_{14}=-(\a_1+2\a_2+3\a_3+5\a_4+4\a_5+3\a_6+2\a_7+\a_8)$,

$-\gamma_{15}=-(2\a_1+2\a_2+3\a_3+4\a_4+4\a_5+3\a_6+2\a_7+\a_8)$,

 It is
easy to check the following fact:

(12) The ten elements $e_{\a_1}$, $e_{\gamma_{11}}, $ $e_{\a_7}$,
$e_{\gamma_2},$ $e_{-\gamma_6}$, $e_{-\gamma_3}$, $e_{\a_2}$,
$e_{\gamma_{12}}$,  $e_{\gamma_{13}}$, $e_{-\gamma_{14}}$ span a
$C_G(s)$-submodule $M_1$ of $\nqs$,

(13) The  eleven elements $e_{\a_3}$, $e_{\gamma_{10}}, $
$e_{-\gamma_5},$ $e_{\a_6}$, $e_{\a_8}$, $e_{-\gamma_8}$,
$e_{-\gamma_{15}},$ $e_{\gamma_1}$, $e_{\a_4}$, $e_{\gamma_9}$,
$e_{-\gamma_4}$ span a $C_G(s)$-submodule $M_2$ of $\nqs$,

(14) $e_{\a_5},e_{-\gamma_7}$ span a $C_G(s)$-submodule $M_3$ of
$\nqs$.

Clearly $\nqs$ is the direct sum of all $M_i$, but it is not easy to
see the number of $C_G(s)$-orbits of $\nqs$ in this case.

If $o(q)=13$,  then $C_G(s)$ is generated by $T, U_{\pm\sigma_i}$,
$i=1,2$, here

$\sigma_1=\a_1+\a_2+2\a_3+3\a_4+3\a_5+2\a_6+\a_7$,

$\sigma_2=\a_1+2\a_2+2\a_3+3\a_4+2\a_5+2\a_6+\a_7$,

$\sigma_3=\a_1+2\a_2+3\a_3+2\a_4+2\a_5+\a_6+2\a_7+\a_8$,

$\sigma_4=\a_1+\a_2+2\a_3+3\a_4+2\a_5+2\a_6+\a_7+\a_8$,

$\sigma_5=\a_1+\a_2+2\a_3+2\a_4+2\a_5+2\a_6+2\a_7+\a_8$,

$\sigma_6=\sigma_1+\sigma_3=2\a_1+3\a_2+4\a_3+6\a_4+5\a_5+3\a_6+2\a_7+\a_8$.

If $\a$ is not simple root, then $\a(s)=q$ if and only if $\a$ is
one of the following roots:

$\tau_1=\a_1+2\a_2+2\a_3+3\a_4+3\a_5+2\a_6+\a_7$,

$-\tau_2=-(\a_1+\a_2+2\a_3+3\a_4+2\a_5+2\a_6+\a_7)$,

$-\tau_3=-(\a_1+2\a_2+2\a_3+3\a_4+2\a_5+\a_6+\a_7)$,

$-\tau_4=-(\a_1+\a_2+2\a_3+3\a_4+2\a_5+\a_6+\a_7+\a_8)$,

$-\tau_5=-(\a_1+\a_2+2\a_3+2\a_4+2\a_5+2\a_6+\a_7+\a_8)$,

$-\tau_6=-(\a_1+\a_2+\a_3+2\a_4+2\a_5+2\a_6+2\a_7+\a_8)$,

$\tau_7=\a_1+2\a_2+2\a_3+3\a_4+2\a_5+2\a_6+\a_7+\a_8$,

$\tau_8=\a_1+\a_2+2\a_3+3\a_4+3\a_5+2\a_6+\a_7+\a_8$,

$\tau_9=\a_1+\a_2+2\a_3+3\a_4+2\a_5+2\a_6+2\a_7+\a_8$,

$-\tau_{10}=-(2\a_1+3\a_2+4\a_3+6\a_4+4\a_5+3\a_6+2\a_7+\a_8)$,

$\tau_{11}=2\a_1+3\a_2+4\a_3+6\a_4+5\a_5+4\a_6+2\a_7+\a_8$.

 It is
easy to check the following facts:

(15) The eleven elements $e_{\a_2}$, $e_{\tau_1}, $ $e_{\a_5}$,
$e_{-\tau_2}$, $e_{-\tau_4}$, $e_{\tau_7}$, $e_{\a_6}$,
$e_{-\tau_3}$, $e_{\tau_{11}}$, $e_{\tau_{10}}$, $e_{\a_8}$ span a
$C_G(s)$-submodule $M_2$ of $\nqs$,

(16) The element $e_{\a_1}$ spans a $C_G(s)$-submodule $M_1$ of
$\nqs$,

(17) The elements $e_{\a_3},\ e_{-\tau_6}$ span a $C_G(s)$-submodule
$M_3$ of $\nqs$.

(18) The elements $e_{\a_4},\ e_{-\tau_5}$, $e_{\tau_9}$, $e_{\a_7}$
span a $C_G(s)$-submodule $M_4$ of $\nqs$.

Clearly $\nqs$ is the direct sum of all $M_i$. But it is not easy to
see the number of $C_G(s)$-orbits of $\nqs$ in this case.

If $o(q)=16$, then $C_G(s)$ is generated by $T, U_{\pm\xi_i}$,
$i=1,2,3,4$, here

$\xi_1=\a_1+2\a_2+3\a_3+4\a_4+3\a_5+2\a_6+\a_7$,

$\xi_2=\a_1+\a_2+2\a_3+3\a_4+3\a_5+3\a_6+2\a_7+\a_8$,

$\xi_3=\a_1+2\a_2+2\a_3+4\a_4+3\a_5+2\a_6+\a_7+\a_8$,

$\xi_4=\a_1+2\a_2+2\a_3+3\a_4+3\a_5+2\a_6+2\a_7+\a_8$.

If $\a$ is not simple root, then $\a(s)=q$ if and only if $\a$ is
one of the following roots:

$\eta_1=2\a_1+2\a_2+3\a_3+4\a_4+3\a_5+2\a_6+\a_7$

$-\eta_2=-(\a_1+2\a_2+2\a_3+4\a_4+3\a_5+2\a_6+\a_7)$,

$-\eta_3=-(\a_1+2\a_2+2\a_3+3\a_4+3\a_5+2\a_6+\a_7+\a_8)$,

$-\eta_4=-(\a_1+2\a_2+2\a_3+3\a_4+2\a_5+2\a_6+2\a_7+\a_8)$,

$-\eta_5=-(\a_1+\a_2+2\a_3+3\a_4+3\a_5+2\a_6+2\a_7+\a_8)$,

$\eta_6=\a_1+2\a_2+2\a_3+3\a_4+3\a_5+3\a_6+2\a_7+\a_8$,

$\eta_7=\a_1+2\a_2+2\a_3+4\a_4+3\a_5+2\a_6+2\a_7+\a_8$,

$\eta_8=\a_1+2\a_2+3\a_3+4\a_4+3\a_5+2\a_6+\a_7+\a_8$.

 It is
easy to check the following facts:

(19) The  elements $e_{\a_1}$, $e_{\eta_1}$ span a
$C_G(s)$-submodule $M_1$ of $\nqs$,

(20) The elements $e_{\a_2}$, $e_{\a_6}$, $e_{\eta_6}$,
$e_{-\eta_5}$ span a $C_G(s)$-submodule $M_2$ of $\nqs$,

(21) The elements $e_{\a_3}$, $e_{\a_8}$, $e_{\eta_8}$,
$e_{-\eta_2}$ span a $C_G(s)$-submodule $M_3$ of $\nqs$,

(22) The elements $e_{\a_4}$, $e_{\a_7}$, $e_{\eta_7}$,
$e_{-\eta_3}$ span a $C_G(s)$-submodule $M_4$ of $\nqs$,

(23) The elements $e_{\a_5},\ e_{-\eta_4}$  span a
$C_G(s)$-submodule $M_5$ of $\nqs$,

Clearly $\nqs$ is the direct sum of all $M_i$. Note that $U_{\pm
\xi_i}$ act on $M_1+M_3$ trivially for $i=2,4$ and $U_{\pm\xi_3}$
also act on $M_1$ trivially. By direct computation we see that the
$C_G(s)$-submodule $M_1+M_3$ has eight $C_G(s)$-orbits,
representatives of the orbits can be chosen as follows: $0,$
$e_{\a_1}$, $e_{\a_3}$, $e_{\a_3}+e_{\a_8}$, $e_{\a_1}+e_{\a_3}$,
$e_{\a_1}+e_{\a_3}+e_{\a_8}$, $e_{\a_1}+e_{-\eta_2}$,
$e_{\a_1}+e_{\a_3}+e_{-\eta_2}$. Similarly, one can see that
$M_4+M_5$ has also eight $C_G(s)$-orbits, and $M_2$ has three
$C_G(s)$-orbits. Thus  the number of $C_G(s)$-orbits in
$\calN_{q,s}$ is not less than $8\times 8\times 3=192$.

If $o(q)=17$, then $C_G(s)$ is generated by $T, U_{\pm\eta_i}$,
$i=1,6,7,8$, here

$\eta_1=2\a_1+2\a_2+3\a_3+4\a_4+3\a_5+2\a_6+\a_7$

$\eta_6=\a_1+2\a_2+2\a_3+3\a_4+3\a_5+3\a_6+2\a_7+\a_8$,

$\eta_7=\a_1+2\a_2+2\a_3+4\a_4+3\a_5+2\a_6+2\a_7+\a_8$,

$\eta_8=\a_1+2\a_2+3\a_3+4\a_4+3\a_5+2\a_6+\a_7+\a_8$.

If $\a$ is not simple root, then $\a(s)=q$ if and only if $\a$ is
one of the following roots:

$-\xi_1=-(\a_1+2\a_2+3\a_3+4\a_4+3\a_5+2\a_6+\a_7)$,

$-\xi_2=-(\a_1+\a_2+2\a_3+3\a_4+3\a_5+3\a_6+2\a_7+\a_8)$,

$-\xi_3=-(\a_1+2\a_2+2\a_3+4\a_4+3\a_5+2\a_6+\a_7+\a_8)$,

$-\xi_4=-(\a_1+2\a_2+2\a_3+3\a_4+3\a_5+2\a_6+2\a_7+\a_8)$,

$\xi_5=2\a_1+2\a_2+3\a_3+4\a_4+3\a_5+2\a_6+\a_7+\a_8$,

$\xi_6=\a_1+2\a_2+3\a_3+4\a_4+3\a_5+2\a_6+2\a_7+\a_8$,

$\xi_7=\a_1+2\a_2+2\a_3+4\a_4+3\a_5+3\a_6+2\a_7+\a_8$.

 It is
easy to check the following facts:

(24) The  elements $e_{\a_2}$, $e_{-\xi_2}$ span a
$C_G(s)$-submodule $M_2$ of $\nqs$,

(25) The element $e_{\a_5}$ spans a $C_G(s)$-submodule $M_5$ of
$\nqs$,

(26) The elements $e_{\a_1},\ e_{-\xi_1}$, $e_{\xi_5}$, $e_{\a_8}$
span a $C_G(s)$-submodule $M_1$ of $\nqs$,

(23) The elements $e_{\a_3},\ e_{-\xi_3}$, $e_{\xi_6}$, $e_{\a_7}$
span a $C_G(s)$-submodule $M_3$ of $\nqs$,

(24) The elements $e_{\a_4},\ e_{-\xi_4}$, $e_{\xi_7}$, $e_{\a_6}$
span a $C_G(s)$-submodule $M_4$ of $\nqs$.

Clearly $\nqs$ is the direct sum of all $M_i$. Noting that $M_i$ has
three $C_G(s)$-orbits for $i=4,5$,  $M_5$ has two $C_G(s)$-orbits,
$M_2+M_4$ has eight $C_G(s)$-orbits, we see that the number of
$C_G(s)$-orbits in $\calN_{q,s}$ is not less than $3\times3\times
8\times2=144$.

Completely similar to the above discussions, one see that the
$C_G(s)$-orbits in $\nqs$ is not less than 144 if
$o(q)=$19,21,22,23,25,26,27,28,29.

\medskip

\centerline{\bf Proof of Theorem 1.2.}

\medskip

In this section we prove  Theorem 1.2, the main result  of this
article.

First assume that $\sum_{w\in W_0}q^{l(w)}\ne 0$.

 Assume that $\psi: \hq\to\cw=H_1$ is an
isomorphism of $\bbC$-algebras. Then $\phi$ induces an isomorphism
from $Z(\hq)\to Z(H_1)$. Let $s$ be a semisimple element in $T$.
Composing the homomorphism $\phi_{1,s}:Z(\cw)\to \bbC$ with $\psi$,
we get a homomorphism $\phi_{1,s}\psi:Z(H_q)\to\bbC$. So there
exists a semisimple element $t$ in $T$ such that
$\phi_{1,s}\psi=\phi_{q,t}$. This induces an isomorphism $H_{q,t}\to
H_{1,s}$.

Assume further that $H_{1,s}$ has one dimensional representation,
then $H_{q,t}$ has also one dimensional representation. So the
number of the central characters $\phi_{q,t}$ $(t\in T)$ which admit
one-dimensional representations is equal to the number of the
central characters $\phi_{1,s}$ $(s\in T)$ which admit
one-dimensional representations. According to Corollary 2.2 and
Corollary 3.3, this is not true if $R$ is not simply laced. So the
theorem is true if $\sum_{w\in W_0}q^{l(w)}\ne 0$ and $R$ is not
simply laced.

\def\ns{|\text{Irr}H_{1,s}|}
\def\nt{|\text{Irr}H_{q,t}|}
\def\nw{|\text{Irr}W_0|}

Now assume that $R$ is simply laced. Let $s,t$ be as above. By 2.3
(a), the number $\ns$ of isomorphism classes of irreducible
representations of $H_{1,s}$ is equal to the number $\nw$ of
isomorphism classes of irreducible complex representations of the
Weyl group $W_0$. We will show that $\nw$ is not equal to  the
number $|\text{Irr}H_{q,t}|$ of isomorphism classes of irreducible
representations of $H_{q,t}$. We do this case by case.

%, $P_2=[\frac k2]+1$ (the integer part of $\frac{k+2}2$)

 Type $A_n$ ($n\ge 2)$. It is known that $\nw$ is the number $p(n+1)$ of
 partitions of $n+1$. Since $\sum_{w\in W_0}q^{l(w)}\ne 0$, the order of $q$ is greater than $n$. By Lemma 3.4, we see
 that $\nt=2^n$. It suffices to show that $\nt=2^n>p(n+1)=\nw$ if $n\ge 2$.

 One can use Euler's Pentagonal Theorem to prove  $2^n>p(n+1)$ for $n\ge 2$. Set $p(0)=1$ and $p(m)=0$ if $m<0$.
 The theorem says
 $$p(n)=\sum_{k=1}^\infty (-1)^{k+1}(p(n-\frac{k(3k-1)}2)+(p(n-\frac{k(3k+1)}2)),$$ which implies that
 $p(n)\le p(n-1)+p(n-2)$ for $n\ge 2$. One  can  also use induction on $n$ to prove
 the result directly. When $n=2,3,4,5$,  $p(n+1)$ is $3, 5,7,11$ respectively. The
 result is true in these cases. Now assume that $k\ge 5$ and the result is true
 for integer $n$ satisfying $2\le n\le k$. For $1\le i\le k+2$, let $P_i$ be the number of partitions
 $a_1\ge a_2\ge \cdots \ge a_m$ of $k+2$
 such that $a_1=i$. Then $P_{k+2}=P_{k+1}=P_1=1,$ $P_{k}=2,$
 $P_{k-1}=3$. By induction hypothesis, $2^{k+1-i}>P_i$ whenever
 $2\le i\le k-2$. Hence, $$\begin{array}{rl}
 p(k+2)&=P_{k+2}+P_{k+1}+P_{k}+P_{k-1}+\cdots+P_2+P_1\\
 &=8+P_{k-2}+P_{k-3}+\cdots+P_3\\
 &<8+2^3+2^4+\cdots +2^{k-1}\\
 &=2^k<2^{k+1}.
\end{array}$$

\def\cp{\mathcal P}
\def\cpn{\cp(n)}
\def\cd{\mathcal D}
\def\cdn{\cd(n)}

 Type $D_n$ $(n\ge 4)$. It is known that $\nw$ is the  number $\mathcal {P}(n)$
 of unordered pairs $(\xi,\eta)$ of
partitions $\xi,\eta$ with $|\xi|+|\eta|=n$, and any pairs
$(\xi,\eta)$ with $\xi=\eta$ and $|\xi|+|\eta|=n$ is counted twice,
see 2.3. By 1.6(b) and the discussion in 3.5, we have $\nt\ge
2^4\cdot 3^{\frac{n-4}2}$ if $n\ge 4$ is even and $\nt\ge 2^5\cdot
3^{\frac{n-5}2}$ if $n\ge 5$ is odd. Define $\cdn$ $(n\ge 4)$ to be
$ 2^4\cdot 3^{\frac{n-4}2}$ if $n$ is even and to be $ 2^5\cdot
3^{\frac{n-5}2}$ if $n$ is odd. It suffices to show that $ \cdn>
\cpn $ if $n\ge 4$.

We first show that

\noindent ($*$) $p(n)\le 2p(n-2)$ if $n\ge 8$.

Let $Q_n$ be the set of all partitions of $n$. For $1\le k\le n$,
let $Q_{n,k}$ be the set of partitions of $n$ with smallest term
 $k$ and $p(n,k)$ be the cardinality of $Q_{n,k}$. Then $p(n,1)=p(n-1)$ and $p(n)$ is the sum of all $p(n,k)$.
 So
$p(n)-p(n-1)$ is the sum of all $p(n,k)$ with $k\ge 2$. Since
$p(n)-p(n-2)=p(n)-p(n-1)+p(n-1)-p(n-2)$, we get
$$p(n)-p(n-2)=\sum_{k=2}^np(n,k)+\sum_{k=2}^{n-1}p(n-1,k).$$

We define an injection
$$\tau:(\bigcup_{k=2}^nQ_{n,k})\cup(\bigcup_{k=2}^{n-1}Q_{n-1,k})\rightarrow\bigcup_{k=1}^{n-2}Q_{n-2,k}=Q_n.$$

For a partition $a_1\ge a_2\ge a_3\cdots \ge a_r>0$ of $n$, we shall
simply write it as $a_1a_2a_3\cdots a_r$.  For $A=a_1a_2\cdots
a_{r-1}a_r\in Q_{n-1,k}$ ($k\ge 2)$, define $\tau(A)=a_1a_2\cdots
a_{r-1}(a_r-1)$. Note that $a_{r-1}>a_r-1\ge 1$.

For $A=a_1a_2\cdots a_{r-1}a_r\in Q_{n,k}$ ($k\ge 4)$, define
$\tau(A)=a_1a_2\cdots a_{r-1}1\cdots   1$, where 1 appears
$a_r-2=k-2$ times. Note that $a_{r-1}-(a_r-2)\ge 2$.

For $A=a_1a_2\cdots a_{r-2}a_{r-1}3\in Q_{n,3}$, set
$\tau(A)=a_1a_2\cdots a_{r-2}21\cdots 1$, where 1 appears
$a_{r-1}-1$ times. Note that $a_{r-1}-1\ge 2$ and $a_{r-2}\ge 3$.

For $r\ge 3$ and $A=a_1a_2\cdots a_{r-3}a_{r-2}a_{r-1}2\in Q_{n,2}$
, set $\tau(A)=a_1a_2\cdots a_{r-3}a_{r-2}a_{r-1}$ if
$a_{r-2}=a_{r-1}$
 and set  $\tau(A)=a_1a_2\cdots a_{r-3}(a_{r-1}+1)1\cdots  1$ if $a_{r-2}>a_{r-1}$, where 1 appears
$a_{r-2}-1$ times. Note that if $a_{r-2}>a_{r-1}$, then
$a_{r-2}-1\ge 2$ and $a_{r-1}+1-(a_{r-2}-1)\le
a_{r-2}-(a_{r-2}-1)\le 1$.

For $A=a_12\in Q_{n,2}$, let $\tau(A)=(a_1-6)2211$. Since $n\ge 8$,
$\tau(A)$ is well defined.

It is easy to check that $\tau$ is injective.  Hence we have
$p(n)-p(n-2)\le p(n-2)$ if $n\ge 8$.

Now we show that $\cdn>\cpn$. We use induction on $n$. When
$n=4,5,6,7,8,9,10,11,12$, $\cpn$ is 13, 18, 37, 55, 100, 150, 251,
376, 599, respectively, and $\cdn$ is 16, 32, 48, 96, 144, 288, 432,
864, 1296, respectively. So the result is true in these cases.

Since $\mathcal D(n+2)=3\cdn$, it suffices to show that $\mathcal
P_{n+2}<3\cpn$ for $n\ge 11$.

\def\cpnt{\mathcal P(n+2)}

Assume that $n=2k\ge 12$. Then
$$\cpn=\sum_{i=0}^{k-1}p(n-i)p(i)+\frac {p(k)(p(k)+3)}2,$$
$$\cpnt=\sum_{i=0}^{k}p(n+2-i)p(i)+\frac {p(k+1)(p(k+1)+3)}2.$$

It is well known that $p(k+1)\le \frac12(p(k+2)+p(k))$. Since
$k+2\ge 8$, by ($*$), we get  $p(k+1)\le \frac32p(k)$. Again by
$(*)$,  for $0\le i\le k-1$, we have $p(n+2-i)\le 2 p(n-i)$ since
$n+2-i=2k+2-i>k+3> 8$. Thus
$$3\cpn-\cpnt\ge \sum_{i=0}^{k-1}p(n-i)p(i)+\frac 38p(k)^2+\frac94p(k)-p(k+2)p(k).$$

Euler's Pentagonal Theorem implies that $p(k)\le p(k-1)+p(k-2)$.
Hence
$$3\cpn-\cpnt\ge \sum_{i=0}^{k-3}p(n-i)p(i)+p(k+1)p(k-1)+\frac
38p(k)^2+\frac94p(k)-p(k+2)p(k-1).$$ But $p(k+2)\le\frac32p(k+1)$,
so $$3\cpn-\cpnt\ge \sum_{i=0}^{k-3}p(n-i)p(i)+\frac
38p(k)^2+\frac94p(k)-\frac12p(k+1)p(k-1).$$ When $k-1\ge 8$, we have
$p(k-1)\le 2p(k-3)$, so $p(k+3)p(k-3)>\frac12p(k+1)p(k-1)$ if $k\ge
9$. If $6\le k<9$, one may check directly that
$p(k+3)p(k-3)>\frac12p(k+1)p(k-1)$. So $3\cpn>3\cpnt$ if $n=2k\ge
12$.

If $n=2k+1\ge 11$, the argument for $3\cpn>\cpnt$ is similar (and
simpler). The proof for type $D_n$ is complete.

Type $E_6,E_7,E_8$. The number $\nw$ is 25, 60, 112, respectively,
see 2.3. Using 1.5 (b), 2.3, Lemma 3.4 and the discussion in 3.5, we
see that $\nt>\nw$ if $o(q)\ne 11, 13$ when $R$ is type $E_8$.

In the rest part of the proof we assume that $\sum_{w\in
W_0}q^{l(w)}= 0$. Then any simple quotient module of certain
standard modules of $H_q$ is a simple constituent of some other
standard modules, see [X1, Theorem 7.8, p.83]. This indicates that
the standard modules of $H_q$ and the standard modules of $H_1$
behave differently,  so $H_q$ and $H_1$ should not be  isomorphic.
We make this point explicit.

Let $t\in T$ be  such that $H_{q,t}$ has a one-dimensional
representation on which all $T_r$ act by scalar $q$. Assume that
$H_q$ and $H_1$ were isomorphic. Then there exists $s\in T$ such
that $H_{q,t}$ and $H_{1,s}$ are isomorphic. By Lemma 3.1, we may
assume that $\a(t)=q$ for all simple roots $\a$. By Lemma 2.1,
$\a(s)=1$ for all simple roots $\a$.

 Let $D=\sum_{w\in
W_0}T_w\in H_q$ and $D'=\sum_{w\in W_0}(-q)^{-l(w)}T_w\in H_q$. We
also use $D$ and $D'$ for their images in $H_{q,t}$ respectively.
Let $C=\sum_{w\in W_0}w\in H_1=\mathbb C[W]$ and $C'=\sum_{w\in
W_0}(-1)^{-l(w)}w\in H_1$. The images in $H_{1,s}$ of $C$ and $C'$
will again be denoted by $C$ and $C'$ respectively.

\def\pai{\pi}
 According to [X3, Theorem 3.3], $CH_{1,s}C'$ is a one-dimensional two-sided ideal of
 $H_{1,s}$.
 Let $\pai: H_{1,s}\to H_{q,t}$ be an isomorphism.
Then $\pai(CH_{1,s}C')=\pi(C)H_{q,t}\pi(C')$ is a one-dimensional
two-sided ideal of $H_{q,t}$. Let $v$ be a nonzero element of
$\pai(CH_{1,s}C')$. Then one of the following cases must happen:

(1) $T_wv=vT_w=q^{l(w)}v$ for all $w\in W_0$,

(2) $T_wv=vT_w=(-1)^{l(w)}v$ for all $w\in W_0$,

(3) $T_wv=q^{l(w)}v,\ vT_w=(-1)^{l(w)}v$ for all $w\in W_0$,

(4) $T_wv=(-1)^{l(w)}v,\ vT_w=q^{l(w)}v$ for all $w\in W_0$.

To go further we need the following facts.

(5) Let $h\in H_q$. If $T_rh=-h$ (resp. $hT_r=-h$; $T_rh=qh$;
$hT_r=qh$) for all simple reflection $r$ in $W_0$, the $h\in D'H_q$
(resp. $h\in H_qD'$; $h\in H_qD$; $h\in DH_q$).

We explain the reasons for (5). For $w\in W$, let $D'_w=\sum_{y\le
w}(-q)^{-l(y)}P_{y,w}(q^{-1})T_y$, where $\le $ stands for the
Bruhat order on $W$, and $P_{y,w}$ are the Kazhdan-Lusztig
polynomials. Then the elements $D'_w, w\in W_0$, form a basis of the
subalgebra $H_{q,W_0}$ of $H_q$ generated by all $T_w, w\in W_0$.
Let $w\in W_0$ and $r$ be a simple reflection in $W_0$. Then
$T_rD'_w=-D'_w$ (resp. $D'_wT_r=-D'_w)$ if and only if $rw\le w$
(resp. $wr\le w$). Moreover, if $rw\ge w$, then
$T_rD'_w=q^{\frac32}D'_{rw}+\sum_{y\le w}a_yD'_y$, $a_y\in
\mathbb{C}$. See (2.3.a) and (2.3.c) in [KL1]. Note that the
elements $D'_w\theta_x,\ w\in W_0,\ x\in X$, form a basis of $H_q$.
Using these facts we see easily that $h\in D'H_q$ if $h\in H_q$ and
$T_rh=-h$ for all simple reflections $r$ in $W_0$. The arguments for
other parts of (5) are similar.

The following statement is easy to check.

\def\t{\theta}
\def\T{\Theta}

(6) Let  $r$ be a simple reflection in $W_0$ and $\t=\sum_{x\in
X}a_x\t_x\in \T_q$, where $a_x\in \mathbb{C}$. If $D\t T_r=qD\t$
(resp. $D'\t T_r=-D'\t$), then $a_x=a_{r(x)}$ and $T_r\t=\t T_r$.

As a consequence of (5) and (6), we get (see also [L1, p.213],  and
note also that $DH_qD\subset DZ(H_q)$ by [L1, Proposition 8.6])

(7) Let  $h\in H_q$. If $T_rh=hT_r=qh$ (resp. $T_rh=hT_r=-h$) for
all simple reflections $r$ in $W_0$, then $h$ is in $DZ(H_q)$ (resp.
$D'Z(H_q)$).

The following result is implicit in  [L1, p.213, line -7]).

(8) Assume that $q\ne -1$. Let  $h\in H_q$. If  $T_rh=-h$ and
$hT_r=qh$ for all simple reflections $r$ in $W_0$, then $h$ is in
$D'H_qD$.

\def\st{\stackrel}
\def\sc{\scriptstyle}

Now we argue for (8).  Assume that $h\ne 0$. By (5), $h$ is in
$H_qD$. Let $\Omega=\{w\in W\,|\, wr\ge w$ for all simple
reflections $r$ in $W_0$\}. According to [X3, Lemma 2.2], the
elements $D'_wD,\ w\in\Omega$, form a basis of $H_qD$. Let
$h=\sum_{w\in \Omega}\xi_wD'_wD,\ \xi_w\in\mathbb{C}$. Assume that
$r$ is a simple reflection in $W_0$. Write $(T_r-q)h=\sum_{w\in
\Omega}\eta_wD'_wD,\ \eta_w\in\mathbb{C}$. Note that for any $w\in
W$ we have $(T_r-q)D'_w=-(1+q)D'_w$ if $rw\le w$ and
$(T_r-q)D'_w=q^{\frac32}D'_{rw}+\sum_{y< w}a_yD'_y$, $a_y\in
\mathbb{C}$ if $rw\ge w$. Moreover, $ry\le y$ if $a_y\ne 0$ (see
[KL1, (2.3.a)]. Since  $D'_yD=0$ for any $y\not\in\Omega$, we must
have $\eta_w=0$ if $rw\ge w$. But $(T_r-q)h=-(1+q)h\ne 0$, so we
must have $rw\le w$ if $\xi_w\ne 0$. Therefore,  $rw\le w$ for all
simple reflections $r$ in $W_0$ whenever $\xi_w\ne 0$.  By (5), then
$D'_w\in D'H_q$ if $\xi_w\ne
 0$. So $h$ is $D'H_qD$.

\medskip

There is a unique involutive automorphism $\xi\to \xi^*$ of the
$\mathbb{C}$-algebra $H_q$ such that $T_r^*=-qT_r^{-1}=q-1-T_r\
(r\in S\cap W_0), \ \theta_x^*=\t_{x^{-1}}\ (x\in X)$ [KL2,
2.13(d)]. Noting that $D^*=(-1)^{l(w_0)}D'$, we see that (8) is
equivalent to the following result.

(9) Let  $h\in H_q$. If  $T_rh=qh$ and $hT_r=-h$ for all simple
reflections $r$ in $W_0$, then $h$ is in $DH_qD'$.

By (7), (8) and (9), we must have $v\in \mathbb CD$, or $v\in
\mathbb CD'$, or $v\in DH_{q,t}D'$, or $v\in D'H_{q,t}D$. But
neither $ \mathbb CD$ nor $ \mathbb CD'$ is two-sided ideal of
$H_{q,t}$. So cases (1) and (2) would not occur. Since $\sum_{w\in
W_0}q^{l(w)}= 0$, by [X3, Theorem 3.3], $ DH_{q,t}D'=D'H_{q,t}D=0$.
These contradict that $v\ne 0$. So $H_{q,t}$ and $H_{1,s}$ are not
isomorphic. Therefore, $H_q$ and $H_1$ are not isomorphic.

The theorem is proved.

\bigskip

\noindent{\bf Acknowledgement:} Part of the work was done during the
visit of N. Xi to the Department of Mathematics, Nagoya University,
Xi thanks the department for hospitality and financial support. We
are very grateful to the referee for carefully reading and valuable
comments, which lead significant improvement of the paper.

\bibliographystyle{unsrt}

\end{document}